\documentclass{amsart}
\usepackage{mathrsfs}
\usepackage{amsfonts}
\usepackage{amssymb}
\usepackage{amsmath,amssymb}

\newtheorem{theorem}{Theorem}[section]
\newtheorem{proposition}[theorem]{Propisition}
\newtheorem{lemma}[theorem]{Lemma}
\newtheorem{corollary}[theorem]{Corollary}
\theoremstyle{definition}

\theoremstyle{remark}
\newtheorem{remark}[theorem]{Remark}

\numberwithin{equation}{section}

\begin{document}

\title[ Note on  affine  Gagliardo-Nirenberg inequalities] {Note on  affine  Gagliardo-Nirenberg  inequalities}


\author{Zhichun Zhai}
\address{Department of Mathematics and Statistics, Memorial University of Newfoundland, St. John's, NL A1C 5S7, Canada}
\curraddr{}
 \email{a64zz@mun.ca}
\thanks{Project supported in part  by Natural Science and
Engineering Research Council of Canada.}

\date{}


\begin{abstract} This note  proves  sharp affine
Gagliardo-Nirenberg inequalities which are stronger than all known
sharp Euclidean Gagliardo-Nirenberg inequalities and imply the
affine $L^{p}-$Sobolev inequalities. The logarithmic version of
affine $L^{p}-$Sobolev inequalities  is verified. Moreover,  An
alternative proof of the affine Moser-Trudinger and Morrey-Sobolev
inequalities is given. The main tools are the equimeasurability of
rearrangements and the strengthened version of the classical
P\'{o}lys-Szeg\"{o} principle.

\end{abstract} \maketitle


 \vspace{0.1in}

 \section{Introduction}
In this note,  we prove sharp affine  Gagliardo-Nirenberg
inequalities. These  inequalities  generalize  the sharp affine
$L^{p}-$Sobolev inequalities
\begin{equation}\label{affine sobolev inequality}
 C_{p,n}\|f\|_{L^{\frac{np}{n-p}}(\mathbb{R}^{n})}
 \leq \mathcal{E}_{p}(f)\ \ \ \
\hbox{for}\ \  f\in W^{1,p}(\mathbb{R}^{n}),
 1\leq p<n,
 \end{equation}
established by Lutwak, Yang and Zhang \cite{Lutwak Yang Zhang 1} for
$1<p<n$ and Zhang \cite{G Zhang} for $p=1.$ Here
$W^{1,p}(\mathbb{R}^{n})$
 is the usual Sobolev space  defined as the set of  functions
 $f\in L^{p}(\mathbb{R}^{n})$ with weak derivative
  $\nabla f\in L^{p}(\mathbb{R}^{n}).$
$\mathcal{E}_{p}(f)$ is the $L^{p}$ affine energy of $f$  defined as
\begin{eqnarray*}
 \mathcal{E}_{p}(f)&=&c_{n,p}\left(\int_{\mathcal{S}^{n-1}}\|D_{v}f\|_{L^{p}(\mathbb{R}^{n})}^{-n}dv\right)^{-\frac{1}{n}}\\
&
=&c_{n,p}\left(\int_{\mathcal{S}^{n-1}}\left(\int_{\mathbb{R}^{n}}|v\cdot\nabla
f(x)|^{p}dx\right)^{-n/p}dv\right)^{-\frac{1}{n}}.
 \end{eqnarray*}
The constant
$c_{n,p}=\left(\frac{n\omega_{n}\omega_{p-1}}{2\omega_{n+p-2}}\right)^{1/p}(n\omega_{n})^{1/n}$
with $\omega_{n}$ being the $n-$dimensional volume enclosed by the
unit sphere $\mathcal{S}^{n-1}.$  For each $v\in \mathcal{S}^{n-1},$
$\|D_{v}f\|_{L^{p}(\mathbb{R}^{n})}$ is the $L^{p}(\mathbb{R}^{n})$
norm of the directional derivative $D_{v}f$ of $f$ along $v.$

   Inequality (\ref{affine sobolev inequality}) is stronger  than the
classical  $L^{p}-$Sobolev inequalities
\begin{equation}\label{sobolev inequality}
 C_{p,n}\|f\|_{L^{\frac{np}{n-p}}(\mathbb{R}^{n})}\leq \|\nabla f\|_{L^{p}(\mathbb{R}^{n})}\ \ \ \  \hbox{for} \ f\in W^{1,p}(\mathbb{R}^{n}),
 \ 1\leq p<n,
 \end{equation}
 see Aubin  \cite{Aubin}
 and Talenti   \cite{G Talenti} for $1<p<n,$
  Federer and Fleming  \cite{Federer Fleming} and Maz'ya  \cite{Mazya}
  for $p=1.$ This can be seen from
 \begin{equation}\label{energey inequality}
 \mathcal{E}_{p}(f)\leq \|\nabla
f\|_{L^{p}(\mathbb{R}^{n})}
\end{equation}
for every $f$ with $\nabla f\in L^{p}(\mathbb{R}^{n})$ and $p\geq1,$
see, Lutwak, Yang and Zhang   \cite{Lutwak Yang Zhang 1}.
  It is well known that
(\ref{sobolev inequality}) does not hold for $p=n$ and $p>n.$ The
Moser-Trudinger inequality and Morrey-Sobolev inequality are
counterparts of  (\ref{sobolev inequality}) for $p=n$ and $p>n,$
respectively. The first one, see
 Moser \cite{Moser},  means that there exits $m_{n}=\sup_{\phi}\int_{0}^{\infty}e^{(\phi(t))^{n'}-t}dt$ such that
\begin{equation}\label{Moser-Trudinger inequality}
\frac{1}{|\hbox{sprt}f|}\int_{\mathbb{R}^{n}}e^{(n\omega_{n}^{1/n}|f(x)|/\|\nabla
f\|_{n})^{n'}}dx\leq m_{n}
\end{equation}
 for every $f\in W^{1,n}(\mathbb{R}^{n})$
with $0<|\hbox{spet}f|:=|\{x\in \mathbb{R}^{n}: f(x)\neq0\}|<\infty$
and $n'=\frac{n}{n-1}.$ Moreover, Carleson and Chang in
\cite{Carleson Chang} proved that extremals do existence for
(\ref{Moser-Trudinger inequality}). Here
 $|A|$ is the Lebesgue measure of $A\subset \mathbb{R}^{n}.$ For
$p>n,$ the Morrey-Sobolev inequality states that
\begin{equation}\label{Morrey sobolev inequality}
\|f\|_{L^{\infty}(\mathbb{R}^{n})}\leq
b_{n,p}|\hbox{sprt}f|^{\frac{1}{n}-\frac{1}{p}}\|\nabla
f\|_{L^{p}(\mathbb{R}^{n})}
\end{equation}
for every $f\in W^{1,p}(\mathbb{R}^{n})$ with
$|\hbox{sprt}f|<\infty.$

 As  a variant of  the classical
$L^{p}-$Sobolev inequality (\ref{sobolev inequality}),
 the Euclidean Gagliardo-Nirenberg /Nash's inequality states that
\begin{equation}
\|f\|_{L^{s}(\mathbb{R}^{n})}\leq C_{n,s,p,q}\|\nabla
f\|^{\theta}_{L^{p}(\mathbb{R}^{n})}\|f\|^{1-\theta}_{L^{q}(\mathbb{R}^{n})}
\end{equation}
for $n\geq1,$  suitable constants $ p, q, s$ and $\theta.$  The
Euclidean Gagliardo-Nirenberg /Nash's inequality has been studied
   intensively and  been applied in analysis and partial differential equations.
 See, for example, Nirenberg \cite{Nirenberg}, Gagliardo \cite{Gagliardo},
  Cordero-Erausquin, Nazaret and Villani \cite{Cordero-Erausquin nazaret
 Villani}, Del Pino  and Dolbeault \cite{Del-Pino Dolbeaut}-\cite{Del-pino
Dolbeault 3} , {Del Pino, Dolbeault and Gentil} \cite{Del-pino
Dolbeault Gentil}, Carlen and Loss \cite{Carlen Loss}, Agueh \cite{M
Agueh}-\cite{M Agueh 2}.

Inequalities  (\ref{Moser-Trudinger inequality}) and (\ref{Morrey
sobolev inequality}) were also strengthened by the  affine
Moser-Trudinger inequality and affine Morrey sobolev inequality (see
Cinachi,  Lutwak, Yang and  Zhang  \cite{A Cianchi E Lutwak D Yang G
Zhang}), respectively.  The main aim of this paper is to  establish
the following sharp affine Gagliardo-Nirenberg. Similar sharp affine
Gagliardo-Nirenberg inequality was studied by Lutwak, Yang and Zhang
in \cite{Lutwak yang Zhang 4}   with the restriction
$s=p\frac{q-1}{p-1}.$ In this paper, we will remove this
restriction. Below,  we will  denote $D^{p,q}(\mathbb{R}^{n})$ as
the completion of the space of smooth compactly supported functions
$f$ on $\mathbb{R}^{n}$ for the norm $\|f\|_{p,q}=\|\nabla
f\|_{L^{p}(\mathbb{R}^{n})}+\|f\|_{L^{q}(\mathbb{R}^{n})}.$

\begin{theorem}\label{theorem general sharp affine GN inequality}
Let $n,p,q$ and $s$ be such that
$$1<p<n \ \ \hbox{and}\ \ 1\leq q<s<p^{*}=\frac{np}{n-p}\ \ \hbox{if} \ \ n>1.$$
Then the $L^{p}$ affine Gagliardo-Nirenberg inequality
\begin{equation}\label{affine general ineq opt GN}
\|f\|_{L^{s}(\mathbb{R}^{n})}\leq
K_{opt}(\mathcal{E}_{p}(f))^{\theta}\|f\|_{L^{q}(\mathbb{R}^{n})}^{1-\theta},
\forall f\in D^{p,q}(\mathbb{R}^{n})
\end{equation}
holds with $\theta=\frac{np(s-q)}{s[np-q(n-p)]},$ and the sharp
constant $K_{opt}>0$ is explicitly given by
$$K_{opt}=\left[\frac{C(n,p,q,s)}{E(u_{\infty})}\right]^{\frac{np+ps-nq}{s[np-q(n-p)]}}.$$
Here
$$C(n,p,q,s)=\frac{\alpha+\beta}{(q\alpha)^{\frac{\alpha}{\alpha+\beta}}(p\beta)^{\frac{\beta}{\alpha+\beta}}}, \alpha=np-s(n-p), \beta=n(s-q)$$
$u_{\infty}$ is the minimizer of the variational problem
\begin{equation}\label{minimizing problem}
\inf\left\{E(u)=\frac{1}{p}\int_{\mathbb{R}^{n}}|\nabla
u|^{p}dx+\frac{1}{q}\int_{\mathbb{R}^{n}}|u|^{q}dx: u\in
D^{p,q}(\mathbb{R}^{n}), \|u\|_{L^{s}(\mathbb{R}^{n})}=1\right\}.
\end{equation}
Moreover,
\begin{equation}\label{extremer general in optial affine GN}
 f_{\sigma,x_{0}}=Cu_{\infty}(
A(x-x_{0}))
\end{equation}
 are optimal functions in inequality (\ref{affine general
ineq opt GN}), for arbitrary $C\neq0,$  $x_{0}\in \mathbb{R}^{n}$
and $A\in GL(n).$
\end{theorem}
\begin{remark}(i) For the proof of  existence of a minimizer to problem (\ref{minimizing
problem}), see, for example,  Del-pino Dolbeault \cite{del-pino
Dolbeault 1}.
\\
(ii) Under the assumption of Theorem \ref{theorem general sharp
affine GN inequality}, (\ref{affine general ineq opt GN})  implies
the $L^{p}$ Gagliardo-Nirenberg inequality, see Agueh \cite{M Agueh
1}-\cite{M Agueh 2}
\begin{equation}\label{general ineq opt GN}
\|f\|_{L^{s}(\mathbb{R}^{n})}\leq K_{opt}\|\nabla
f\|_{L^{p}(\mathbb{R}^{n})}^{\theta}\|f\|_{L^{q}(\mathbb{R}^{n})}^{1-\theta},
\forall f\in D^{p,q}(\mathbb{R}^{n}).
\end{equation}
 Moreover, $f_{\sigma,x_{0}}=Cu_{\infty}(\sigma
(x-x_{0}))$ are optimal functions in inequality (\ref{affine general
ineq opt GN}), for arbitrary $C\neq0,$ $\sigma\neq0$ and $x_{0}\in
\mathbb{R}^{n}.$ \\
 (iii)If $q=1$ and $p=s$, from (\ref{affine general
ineq opt GN}), we can get the affine $L^{p}$ Nash's inequality
$$\left(\int_{\mathbb{R}^{n}}|f(x)|^{p}dx\right)^{1+\frac{p}{n(p-1)}}\leq (K_{opt})^{p+\frac{p^{2}}{n(p-1)}}(\mathcal{E}_{p}(f))^{p}
\left(\int_{\mathbb{R}^{n}}|f(x)|dx\right)^{\frac{p^{2}}{n(p-1)}}$$
for $1<p<n$ if $n>1.$
\end{remark}

Theorem \ref{theorem general sharp affine GN inequality} implies the
following sharp affine Gagliardo-Nirenberg inequalities stronger
than the Euclidean ones  in  \cite{del-pino Dolbeault 1}.
\begin{corollary}\label{theorem opt affine GN}
Let $1<p<n,$ $p<q\leq \frac{p(n-1)}{n-p}.$ Then for all $f\in
D^{p,q}(\mathbb{R}^{n}),$  we have
\begin{equation}\label{affine ineq opt GN}
\|f\|_{L^{s}(\mathbb{R}^{n})}\leq
C_{2}(\mathcal{E}_{p}(f))^{\theta}\|f\|_{L^{q}(\mathbb{R}^{n})}^{1-\theta}.
\end{equation}
Here $s=p\frac{q-1}{p-1}$ and
$$\theta=\frac{(q-p)n}{(q-1)(np-(n-p)q)}$$
and with $\delta=np-q(n-p)>0,$ the optimal constant $C_{2}$ takes
the form
$$C_{2}=\left(\frac{q-p}{p\sqrt{\pi}}\right)^{\theta}\left(\frac{pq}{n(q-p)}\right)^{\frac{\theta}{p}}
\left(\frac{\delta}{pq}\right)^{\frac{1}{s}}\left(\frac{\Gamma\left(q\frac{p-1}{q-p}\right)\Gamma\left(\frac{n}{2}+1\right)}
{\Gamma\left(\frac{(p-1)}{p}\frac{\delta}{q-p}\right)\Gamma\left(n\frac{p-1}{p}+1\right)}\right)^{\frac{\theta}{n}}.$$
Equality holds in (\ref{affine ineq opt GN}) if and only if for some
$\alpha\in \mathbb{R},$ $\beta>0,$ $\overline{x}\in \mathbb{R}^{n},$
\begin{equation}\label{extremer in optial affine GN}
f(x)=\alpha\left(1+\beta|A(x-\overline{x})|^{\frac{p}{p-1}}\right)^{-\frac{p-1}{q-p}}\
\ \ \ \forall x\in \mathbb{R}^{n}
\end{equation}
with $A\in GL(n).$
\end{corollary}
\begin{remark}
(i) When $q=p\frac{n-1}{n-p},$ $\theta=1$ and $s=\frac{np}{n-p}.$
Thus inequality (\ref{affine ineq opt GN}) implies the sharp  affine
$L^{p}-$Sobolev inequality.\\
 (ii) Inequality (\ref{affine ineq opt GN}) was proved by Lutwak, Yang and  Zhang in \cite{Lutwak yang Zhang 4}
 where the authors applied the optimal $L^{p}$ Sobolev norm problems
 and $L^{p}$ Petty projection inequality (see Gardner \cite{Gardner}, Schneider \cite{Schneider} and Thompson \cite{A.C. Thompson}
 for  $p=1,$ Lutwak, Yang  and Zhang \cite{Lutwak Yang Zhang} for $p>1.$)
  \end{remark}

Similarly, for $q<p<n,$  we can obtain the following resutls.

\begin{corollary}\label{theorem another opt affine GN}
Let $1<p<n,$ $1<q<p.$ Then for all $f\in D_{p,r}(\mathbb{R}^{n}),$
we have
\begin{equation}\label{affine another ineq opt GN}
\|f\|_{L^{q}(\mathbb{R}^{n})}\leq
C_{3}(\mathcal{E}_{p}(f))^{\theta}\|f\|_{L^{r}(\mathbb{R}^{n})}^{1-\theta}.
\end{equation}
Here $r=p\frac{q-1}{p-1}$ and
$$\theta=\frac{(p-q)n}{q(n(n-q)+p(q-1))}$$
and with $\delta=np-q(n-p)>0,$ the optimal constant $C_{3}$ takes
the form
$$C_{3}=\left(\frac{p-q}{p\sqrt{\pi}}\right)^{\theta}\left(\frac{pq}{n(p-q)}\right)^{\frac{\theta}{p}}
\left(\frac{pq}{\delta}\right)^{\frac{1-\theta}{r}}\left(\frac{\Gamma\left(\frac{p-1}{p}\frac{\delta}{p-q}+1\right)\Gamma\left(\frac{n}{2}+1\right)}
{\Gamma\left(q\frac{(p-1)}{p-q}+1\right)\Gamma\left(n\frac{p-1}{p}+1\right)}\right)^{\frac{\theta}{n}}.$$
If $q>2-\frac{1}{p},$ equality holds in (\ref{affine ineq opt GN})
if and only if for some $\alpha\in \mathbb{R},$ $\beta>0,$
$\overline{x}\in \mathbb{R}^{n},$
\begin{equation}\label{extremer in optial affine GN}
f(x)=\alpha\left(1-\beta|A(x-\overline{x})|^{\frac{p}{p-1}}\right)_{+}^{-\frac{p-1}{q-p}}\
\ \ \ \forall x\in \mathbb{R}^{n}
\end{equation}
with $A\in GL(n).$
\end{corollary}

We get the following logarithmic version of (\ref{affine sobolev
inequality}).

\begin{proposition}\label{theorem affine ineq Sobo loga}  For any $f\in W^{1,p}(\mathbb{R}^{n})$ with
 $1<p<n$ and $\int_{\mathbb{R}^{n}}|f(x)|^{p}dx=1,$ we have
\begin{equation}\label{affine ineq Sobo loga}
\int_{\mathbb{R}^{n}}|f(x)|^{p}\log|f(x)|dx\leq\frac{n}{p^{2}}\log\left(C_{4}(\mathcal{E}_{p}(f))^{p}\right).
\end{equation}
 Here the optimal constant $C_{4}$ is defined by
  \begin{equation}\label{constant in ineq Sobo loga}
C_{4}=\frac{p}{n}\left(\frac{p-1}{e}\right)^{p-1}\pi^{-\frac{p}{2}}\left[\frac{\Gamma\left(\frac{n}{2}+1\right)}
{\Gamma\left(n\frac{p-1}{p}+1\right)}\right]^{\frac{p}{n}}.
\end{equation}
 Inequality in (\ref{affine ineq Sobo loga}) is
optimal and equality holds if and only if for some $\sigma>0$ and
$\overline{x}\in \mathbb{R}^{n},$
\begin{equation}\label{extremer in affine Sobo loga}
f(x)=\pi^{\frac{n}{2}}\sigma^{-n\frac{(p-1)}{p}}\frac{\Gamma\left(\frac{n}{2}+1\right)}{\Gamma\left(n\frac{(p-1)}{p}+1\right)}
e^{-\frac{1}{\sigma}|A(x-\overline{x})|^{\frac{p}{p-1}}} \ \ \ \ \ \
\forall x\in \mathbb{R}^{n}
\end{equation}
with $A\in GL(n).$
\end{proposition}

\begin{remark}
 Inequality (\ref{affine ineq Sobo loga}) generalizes  the sharp
Euclidean $L^{p}-$Sobolev logarithmic equality  since
$\mathcal{E}_{p}(f)\leq \|\nabla f\|_{L^{p}(\mathbb{R}^{n})}.$
Meanwhile, it can also been viewed as the limiting case $r=p=q$ of
inequality (\ref{affine ineq opt GN}).  For more details about
Euclidean $L^{p}-$Sobolev logarithmic equality,
 see Weissler \cite{Weissler} and  Groos \cite{Gross},
   Del Pino and  Dolbeault \cite{Del-Pino Dolbeaut}, Gentil \cite{Gentil} and the reference therin.
\end{remark}

 We give an alternative  proof of the affine Moser-Trudinger and Morrey-Sobolev
 inequalities established by
  Cianchi, Lutwak, Yang  and Zhang in \cite{A Cianchi E Lutwak D Yang G
  Zhang}.
  \begin{proposition}\label{proposition moser trudinger}
  Suppose $n>1.$ Then for every $f\in W^{i,n}(\mathbb{R}^{n})$ with
  $0<|\hbox{supp}(f)|<\infty,$
  \begin{equation}\label{sharp affine moser trudinger}
\frac{1}{|\hbox{supp}(f)|}\int_{\hbox{supp}(f)}\exp\left(n\omega_{n}\frac{|f(x)|}{\mathcal{E}_{n}(f)}\right)^{n'}dx\leq
m_{n}
  \end{equation}
with $m_{n}=\sup_{\phi}\int_{0}^{\infty}e^{(\phi(t))^{n'}-t}dt.$ The
constant $n\omega_{n}^{1/n}$ is the best one in the sense that
(\ref{sharp affine moser trudinger}) would fail if
$n\omega_{n}^{1/n}$ is replaced by a larger one.
  \end{proposition}
\begin{proposition}\label{proposition morrey sobolev}
If $p<n,$ then for every $f\in W^{1,p}(\mathbb{R}^{n})$ with
$|\hbox{sprt}(f)|<\infty,$
\begin{equation}\label{affine morrey sobolev}
\|f\|_{L^{\infty}(\mathbb{R}^{n})}\leq
b_{n,p}|\hbox{sprt}(f)|^{\frac{1}{n}-\frac{1}{p}}\mathcal{E}_{p}(f).
\end{equation}
Equality holds in (\ref{affine morrey sobolev}) if and only if
$$f(x)=a\left(1-|A(x-x_{0})|)^{\frac{p-n}{p-1}}\right)_{n}$$
for some $a\in \mathbb{R},$ $x_{0}\in\mathbb{R}^{n},$ and $A\in
GL(n).$ Here $``+"$ denotes the  ``positive part".
\end{proposition}
  Cianchi, Lutwak, Yang  and Zhang,  in \cite{A Cianchi E Lutwak D Yang G
  Zhang}, proved inequality (\ref{sharp affine moser trudinger}) by showing that
$$m_{n}=\sup_{\phi}\frac{1}{a}\int_{0}^{a}\exp(n\omega_{n}^{1/n}\phi(s))^{n'}ds$$
and inequality (\ref{affine morrey sobolev}) by the the strengthened
version of the classical P\'{o}lya-Szeg\"{o} principle,  the local
absolute continuity of the  decreasing rearrangement of $f$ and the
H\"{o}lder inequality. Here, we will prove inequalities (\ref{sharp
affine moser trudinger}) and (\ref{affine morrey sobolev}) directly
by the observation that sphere   rearrangements of functions may
give us better estimates for (affine) Sobolev type inequalities.

The rest of this paper is organized as follows: In Section 2, we
recall some basic properties of rearrangements of  functions
 and the strengthened
version of the classical P\'{o}lya-Szeg\"{o} principle. In Section
3, we prove Propositions   \ref{theorem opt affine GN}-
\ref{proposition morrey sobolev}.

\section{Strengthened Version of the Classical  P\'{o}lya-Szeg\"{o} Principle}

Let $f:\mathbb{R}^{n}\longrightarrow\mathbb{R}$ with
\begin{equation}\label{vanishing at infnity}
|\{x\in\mathbb{R}^{n}: |f(x)|>t\}|<\infty\ \ \  \hbox{for}\ t>0.
\end{equation}
 The distribution function $m_{f}(t)$ of $f$ is defined as
$$m_{f}(t)=|\{x\in \mathbb{R}^{n}:|f(x)|>t\}|, \ \ \ \hbox{for}\
t\geq0.
$$
Functions having the same distribution function are refered to be
equidistributed or equimeasurable. On the other hand,
equidistributed functions are said to be rearrangements of  each
other. The decreasing rearrangement $f^{*}$ of function $f$ is
defined as
$$f^{*}(s)=\sup\{t\geq 0:m_{f}(t)>s\}\ \ \ \ \hbox{for}\ \ s\geq 0.$$
The spherical symmetric rearrangement
$f^{\star}:\mathbb{R}^{n}\longrightarrow[0,\infty]$ is defined as
$$f^{\star}(x)=f^{*}(\omega_{n}|x|^{n})\ \ \ \ \hbox{for}\
 \ x\in \mathbb{R}^{n}.$$
Clearly, $f,f^{*}$ and $f^{\star}$ are equidistributed functions. In
fact, we have
$$m_{f}=m_{f^{*}}=m_{f^{\star}},$$
\begin{equation}\label{spport}
|\hbox{sprt}(f)|=|\hbox{sprt}(f^{*})|=|\hbox{sprt}(f^{\star})|,
\end{equation}
\begin{equation}\label{L infty norm}
\|f\|_{L^{\infty}(\mathbb{R}^{n})}=f^{*}(0)=\|f^{\star}\|_{L^{\infty}},
\end{equation}
and
\begin{equation}\label{general equaimeasurablity}
\int_{\mathbb{R}^{n}}\Phi(|f(x)|)dx=\int_{0}^{\infty}\Phi(f^{*}(s))ds
=\int_{\mathbb{R}^{n}}\Phi(f^{\star}(x))dx
\end{equation}
for every continuous increasing function
$\Phi:[0,\infty)\longrightarrow[0,\infty).$
 Thus, we have
\begin{equation}\label{equaimeasurablity}
\int_{\mathbb{R}^{n}}\|f(x)|^{p}dx=\int_{0}^{\infty}(f^{*}(s))^{p}ds
=\int_{\mathbb{R}^{n}}(f^{\star}(x))^{p}dx
\end{equation}
for every $p\geq 1,$ and so Lebesgue norms will be invariant under
the operations of decreasing rearrangement and of spherically
symmetric rearrangement.

The classical P\'{o}lya-Szeg\"{o} principle means that if $f$ with
(\ref{vanishing at infnity}), is weakly differentiable in
$\mathbb{R}^{n}$ and $|\nabla f|\in L^{p}(\mathbb{R}^{n})$ for
$p\in[1,\infty],$ then $f^{*}$ is locally absolutely continuous in
$(0,+\infty),$ $f^{\star}$ is weakly differentiable in
$\mathbb{R}^{n}$ and
\begin{equation}\label{polya szego}
\|\nabla f^{\star}\|_{L^{p}(\mathbb{R}^{n})}\leq \|\nabla
f\|_{L^{p}(\mathbb{R}^{n})}.
\end{equation}
See, for example,  Kawohl \cite{B Kawohl}-\cite{B Kawohl1}, Sperner
\cite{E Sperner},  Talenti \cite{G Talenti},  Brothers and Ziemer
\cite{JE Brothers WP Ziemer},  Hilden \cite{K Hilden}. Inequality
(\ref{polya szego}) is a powerful tool to many  problems in  physics
and mathematics. On the other hand, several variants of inequality
(\ref{polya szego}) have been established and applied intensively,
see, for example, Kawohl \cite{B Kawohl1}.  Especially, Lutwak, Yang
and Zhang in \cite{Lutwak Yang Zhang},   Cianchi, Lutwak, Yang and
Zhang in \cite{A Cianchi E Lutwak D Yang G Zhang} proved the
following strengthened affine version of inequality (\ref{polya
szego}).

\begin{lemma} \label{lemma affine polya szego}\cite{Lutwak Yang Zhang} \cite{A Cianchi E Lutwak D Yang G Zhang}
Suppose $1<n$ and $1\leq p.$ If $f\in W^{1,p}(\mathbb{R}^{n}),$ then
$f^{\star}\in W^{1,p}(\mathbb{R}^{n}),$
\begin{equation}\label{affine polya szego}
\mathcal{E}_{p}(f^{\star})\leq \mathcal{E}_{p}(f)
\end{equation}
and
\begin{equation}\label{affine equality polya szego}
\|\nabla f^{\star}\|=\mathcal{E}_{p}(f^{\star}).
\end{equation}
\end{lemma}
\begin{remark} \label{remark lemma affine polya szego} 
We can see that both (\ref{affine polya szego}) and (\ref{affine
equality polya szego}) are true for $f\in D^{p,q}(\mathbb{R}^{n}).$
\end{remark}

Inequality (\ref{polya szego}) is particular  significant for the
authors in  \cite{G Zhang}, \cite{Lutwak Yang Zhang 1} and \cite{A
Cianchi E Lutwak D Yang G Zhang} to proved the affine
$L^{p}-$Sobolev, affine Moser-Trudinger and affine Morrey-Sobolev
inequalities. In this note, we will see that inequality (\ref{polya
szego})  implies
 the affine Gagliardo-Nirenberg inequalities.

The proof of Lemma \ref{lemma affine polya szego} depends on $L^{p}$
Brunn-Minkowsi theory of convex bodies (see, for example, Chen
\cite{W Chen}, Chou and Wang \cite{KS Chou XJ Wang},  Hu, Ma  and
Shen \cite{C Hu XN Ma C Shen}, Ludwig \cite{M Ludwig}-\cite{M Ludwig
2}, Lutwak \cite{E Lutwak 1}-\cite{E. Lutwak},  Lutwak and Oliker
\cite{E Lutwak V Oliker}, Lutwak, Yang and Zhang \cite{Lutwak Yang
Zhang}-\cite{Lutwak yang Zhang 4}). In \cite{A Cianchi E Lutwak D
Yang G Zhang}, Lutwak, Yang and Zhang proved Lemma \ref{lemma affine
polya szego}  by applying the similar rearrangement argument used to
prove the Euclidean Sobolev inequality. They solved a family of
$L^{p}$ Minkowski problem (see, Lutwak, Yang and Zhang \cite{Lutwak
yang Zhang 3}) to reduce the estimates for $L^{p}$ gradient
integrals to the estimates for $L^{p}$ mixed volumes of convex
bodies. Thus they can  replace the classical Euclidean isoperimetric
inequality by the affine $L^{p}$ isoperimetric inequality (see,
Lutwak, Yang and Zhang \cite{Lutwak Yang Zhang}). For the details of
 Lemma \ref{lemma affine polya szego}, we refer the
interested  reader to Lutwak, Yang and Zhang \cite[Theorem 2.1]{A
Cianchi E Lutwak D Yang G Zhang}.

\section{Proof of the Main Results}
\subsection{Proof of Theorem \ref{theorem general sharp affine GN inequality}}

The symmetrization inequality (\ref{affine polya szego}) and
inequality (\ref{equaimeasurablity}) are crucial for the proof of
Theorem \ref{theorem general sharp affine GN inequality}.

For any $f\in D^{p,q}(\mathbb{R}^{n}),$
 inequality
(\ref{equaimeasurablity}) implies that
\begin{eqnarray*}
\|f\|_{L^{q}(\mathbb{R}^{n})}
=\|f^{\star}\|_{L^{q}(\mathbb{R}^{n})}.
\end{eqnarray*}
The classical P\'{o}lya-Szeg\"{o} principle  and inequality
(\ref{equaimeasurablity}) tell us that $f^{\star}\in
D^{p,q}(\mathbb{R}^{n}).$ Thus, we can apply the sharp Euclidean
Gagliardo-Nirenberg inequality (\ref{general ineq opt GN})(see
\cite[Theorem 2.1]{M Agueh 2}) to $f^{\star}$ and get
\begin{eqnarray*}\|f^{\star}\|_{L^{s}(\mathbb{R}^{n})}\|f^{\star}\|
_{L^{q}(\mathbb{R}^{n})}^{\theta-1}\leq C_{2}\|\nabla
f^{\star}\|^{\theta}_{L^{p}(\mathbb{R}^{n})}.
 \end{eqnarray*}
Lemma \ref{lemma affine polya szego} and  Remark \ref{remark lemma
affine polya szego} imply
\begin{eqnarray*}
\|f\|_{L^{s}(\mathbb{R}^{n})}\|f\|_{L^{q}(\mathbb{R}^{n})}^{\theta-1}
&=&\|f^{\star}\|_{L^{s}(\mathbb{R}^{n})}\|f^{\star}\|_{L^{q}(\mathbb{R}^{n})}^{\theta-1}\\
&\leq& K_{opt}\|\nabla f^{\star}\|_{L^{p}(\mathbb{R}^{n})}^{\theta}\\
&=& K_{opt}(\mathcal{E}_{p}(f^{\star}))^{\theta}\\
&\leq&K_{opt}(\mathcal{E}_{p}(f))^{\theta}.
\end{eqnarray*}
Thus, we get
$$\|f\|_{L^{s}(\mathbb{R}^{n})}\leq K_{opt}(\mathcal{E}_{p}(f))^{\theta}\|f\|_{L^{q}(\mathbb{R}^{n})}^{1-\theta}.$$
On the other hand, since
\begin{equation}\label{inequality in proof of theorem 1}
\|f\|_{L^{s}(\mathbb{R}^{n})}\leq
K_{opt}(\mathcal{E}_{p}(f))^{\theta}\|f\|_{L^{q}(\mathbb{R}^{n})}^{1-\theta}\leq
 K_{opt}\|\nabla f\|_{L^{p}(\mathbb{R}^{n})}^{\theta}\|f\|_{L^{q}(\mathbb{R}^{n})}^{1-\theta},
 \end{equation}
 the extremal  for
 sharp Gagliardo-Nirenberg inequality (\ref{general ineq opt GN})
 is an extremal of (\ref{inequality in proof of theorem 1}).
 It is easy to see that inequality (\ref{affine ineq opt GN}) is
 an affine inequality, thus composing the extremal functions of inequality (\ref{general ineq opt GN})
  with an element from $GL(n)$ will also give an extremal for the
 affine Gagliardo-Nirenberg inequality (\ref{affine general  ineq opt GN}).
 Thus, the function given by  (\ref{extremer general in optial affine GN}) is the
 extremal of inequality (\ref{affine general ineq opt GN}).

\subsection{Proof of Proposition \ref{theorem affine ineq Sobo loga}}
We combine the symmetrization inequality (\ref{affine polya szego})
and inequality (\ref{general equaimeasurablity}) to prove  Theorem
\ref{theorem affine ineq Sobo loga}. Since $G(t)=t^{p}\log
t:[0,\infty)\longrightarrow[0,\infty)$ is continuous increasing,
inequality (\ref{general equaimeasurablity}) verifies
$$\int_{\mathbb{R}^{n}}|f(x)|^{p}\log|f(x)|dx=\int_{\mathbb{R}^{n}}|f^{\star}(x)|^{p}\log|f^{\star}(x)|dx.$$
Lemma \ref{lemma affine polya szego} verifies   $f^{\star}\in
W^{1,p}(\mathbb{R}^{n})$ and inequality (\ref{equaimeasurablity})
implies
$\|f^{\star}\|_{L^{p}(\mathbb{R}^{n})}=\|f\|_{L^{p}(\mathbb{R}^{n})}=1.$
Thus, we can apply   the sharp Euclidean $L^{p}-$Sobolev logarithmic
inequality (see Del Pino Dolbeaut \cite[Theorem 1.1]{Del-Pino
Dolbeaut}) to $f^{\star}$ and obtain
$$\int_{\mathbb{R}^{n}}|f^{\star}(x)|^{p}\log|f^{\star}(x)|dx
\leq \frac{n}{p^{2}}\log\left(C_{4}\|\nabla
f^{\star}\|^{p}_{L^{p}(\mathbb{R}^{n})}\right).$$ Similar to the
proof of Theorem \ref{theorem opt affine GN}, Lemma \ref{lemma
affine polya szego} gives us
\begin{eqnarray*}
\int_{\mathbb{R}^{n}}|f(x)|^{p}\log|f(x)|dx
&=&\int_{\mathbb{R}^{n}}|f^{\star}(x)|^{p}\log|f^{\star}(x)|dx\\
&\leq& \frac{n}{p^{2}}\log\left(C_{4}\|\nabla
f^{\star}\|^{p}_{L^{p}(\mathbb{R}^{n})}\right)\\
&=&
\frac{n}{p^{2}}\log\left(C_{4}(\mathcal{E}_{p}(f^{\star}))^{p}\right)\\
&\leq&
\frac{n}{p^{2}}\log\left(C_{4}(\mathcal{E}_{p}(f))^{p}\right).
\end{eqnarray*}
Thus, $$\int_{\mathbb{R}^{n}}|f(x)|^{p}\log|f(x)|dx \leq
\frac{n}{p^{2}}\log\left(C_{4}(\mathcal{E}_{p}(f))^{p}\right).$$ The
 function given
by (\ref{extremer in affine Sobo loga}) is an extremal function
inequality (\ref{affine ineq Sobo loga}) since it is also an
extremal function of the sharp Euclidean $L^{p}-$Sobolev inequality
and $$\int_{\mathbb{R}^{n}}|f(x)|^{p}\log|f(x)|dx \leq
\frac{n}{p^{2}}\log\left(C_{4}(\mathcal{E}_{p}(f))^{p}\right)
\leq\frac{n}{p^{2}}\log\left(C_{4}(\|\nabla
f\|_{L^{p}(\mathbb{R}^{n})})^{p}\right).$$

\subsection{Proof of Proposition \ref{proposition moser trudinger}}
Under the assumption of Proposition \ref{proposition moser
trudinger}, the sharp Moser-Trudinger inequality
(\ref{Moser-Trudinger inequality}) holds.  It follows from
(\ref{energey inequality}) and Lemma \ref{affine polya szego} that
$$\|\nabla
f^{\star}\|_{L^{n}(\mathbb{R}^{n})}=\mathcal{E}_{n}(f^{\star})\leq
\mathcal{E}_{n}(f)\leq \|\nabla f\|_{L^{n}(\mathbb{R}^{n})}.
$$
Then we get
\begin{eqnarray*}
&&\frac{1}{|\hbox{supp}(f)|}\int_{|\hbox{supp}(f)|}\exp\left(n\omega_{n}^{1/n}\frac{|f(x)|}{\mathcal{E}_{n}(f)}\right)^{n'}dx\\
&\leq&\frac{1}{|\hbox{supp}(f)|}\int_{|\hbox{supp}(f)|}\exp\left(n\omega_{n}^{1/n}\frac{|f(x)|}{\mathcal{E}_{n}(f^{\star})}\right)^{n'}dx\\
&=&\frac{1}{|\hbox{supp}(f^{\star})|}\int_{|\hbox{supp}(f^{\star})|}\exp\left(n\omega_{n}^{1/n}\frac{|f^{\star}(x)|}{\|\nabla f^{\star}\|_{L^{n}(\mathbb{R}^{n})}}\right)^{n'}dx\\
&\leq&m_{n}
\end{eqnarray*}
with  the last  inequality using  (\ref{Moser-Trudinger
inequality}).
 This implies that (\ref{sharp affine moser
trudinger}) holds. Since \begin{eqnarray*}
&&\frac{1}{|\hbox{supp}(f)|}\int_{|\hbox{supp}(f)|}\exp\left(n\omega_{n}^{1/n}\frac{|f(x)|}{\|\nabla
f\|_{L^{n}(\mathbb{R}^{n})}}\right)^{n'}dx\\
&\leq&
\frac{1}{|\hbox{supp}(f)|}\int_{|\hbox{supp}(f)|}\exp\left(n\omega_{n}^{1/n}\frac{|f(x)|}{\mathcal{E}_{n}(f)}\right)^{n'}dx\leq
m_{n}
\end{eqnarray*}
and  extremal functions for (\ref{Moser-Trudinger inequality})
exist, we see that $f(Ax)$ is an extremal function of (\ref{sharp
affine moser trudinger}) for every  extremal function $f$ for
(\ref{Moser-Trudinger inequality}) and $A\in GL(n).$ On the other
hand, to see the sharpness of $n\omega_{n}^{1/n},$ we assume that
(\ref{sharp affine moser trudinger}) is true for some $\beta>
n\omega_{n}^{1/n}$ and any $f\in W^{1,n}(\mathbb{R}^{n})$ with
$0<|\hbox{supp}(f)|<\infty.$ Then we have $f^{\star}\in
W^{1,n}(\mathbb{R}^{n})$ and
\begin{eqnarray*}
&&\frac{1}{|\hbox{supp}(f^{\star})|}\int_{|\hbox{supp}(f^{\star})|}\exp\left(\beta\frac{|f^{\star}(x)|}{\mathcal{E}_{n}(f^{\star})}\right)^{n'}dx\\
&=&\frac{1}{|\hbox{supp}(f^{\star})|}\int_{|\hbox{supp}(f^{\star})|}\exp\left(\beta\frac{|f^{\star}(x)|}{\|\nabla f^{\star}\|_{L^{n}(\mathbb{R}^{n})}}\right)^{n'}dx\\
&\leq&m_{n}.
\end{eqnarray*}
The last inequality contradicts with the sharpness of
$n\omega_{n}^{1/n}$ in (\ref{Moser-Trudinger inequality}). This
finishes the proof of Proposition \ref{proposition moser trudinger}.

\subsection{Proof of Proposition \ref{proposition morrey sobolev}}
Assume that $f\in W^{1,p}$ with $|\hbox{sprt}{f}|<\infty.$ Then,
from the classical P\'{o}lya-Szeg\"{o} principle, we know that
$f^{\star}\in W^{1,p}(\mathbb{R}^{n}).$ On the other hand, equality
(\ref{spport}) implies that $|\hbox{sprt}(f^{\star})|<\infty.$ Thus,
for $f^{\star},$  we can apply the classical Morrey-Sobolev
inequality and get
$$\|f^{\star}\|_{L^{\infty}(\mathbb{R}^{n})}\leq b_{n,p}
|\hbox{sprt}(f^{\star})|^{\frac{1}{n}-\frac{1}{p}}\|\nabla
f\|_{L^{p}(\mathbb{R}^{n})}.$$
 Equality (\ref{L infty norm}) and Lemma \ref{lemma affine polya szego} imply
 that
\begin{eqnarray*}
\|f\|_{L^{\infty}(\mathbb{R}^{n})}=\|f^{\star}\|_{L^{\infty}}&\leq&
b_{n,p} |\hbox{sprt}(f^{\star})|^{\frac{1}{n}-\frac{1}{p}}\|\nabla
f^{\star}\|_{L^{p}(\mathbb{R}^{n})}\\
&=& b_{n,p}
|\hbox{sprt}(f^{\star})|^{\frac{1}{n}-\frac{1}{p}}\mathcal{E}_{p}(f^{\star})\\
&\leq& b_{n,p}
|\hbox{sprt}(f^{\star})|^{\frac{1}{n}-\frac{1}{p}}\mathcal{E}_{p}(f).
\end{eqnarray*}
The verifying of extremal functions is obviously since the affine
invariance of (\ref{affine morrey sobolev}).

 \vspace{0.1in} \noindent
 {\bf{Acknowledgements.}} The author
  would like to thank  Professor Jie Xiao
 for all kind  encouragement.


\begin{thebibliography}{10}
\bibitem{M Agueh}
M. Agueh, \textit{Sharp Gagliardo-Nirenberg inequalities and Mass
transport theory,} {Journal of Dynamics and Differential Equations}
\textbf{18}(4) (2006), 1069-1093.

\bibitem{M Agueh 1}
 M. Agueh, \textit{Gagliardo-Nirenberg
inequalities involving the gradient $L^{2}-$norm,} {C. R. Math.
Acad. Sci. Paris, Ser. I} \textbf{346} (2008), 757-762.

\bibitem{M Agueh 2}
 M. Agueh,
\textit{Sharp Gagliardo-Nirenberg Inequalities via p-Laplacian Type
Equations,} {Nonlinear differ. equ. appl.} \textbf{15} (2008),
457-472.

\bibitem{Aubin}
T. Aubin, \textit{Probl\`{e}mes isop\'{e}rimetriques et espaces de
Sobolev,} {J. Differ. Geom.} \textbf{11} (1976), 573¨C598.

\bibitem{JE Brothers WP Ziemer}
J.E. Brothers, W.P. Ziemer, \textit{Minimal rearrangements of
Sobolev functions,} {J. Reine Angew. Math} \textbf{384} (1988),
153-179.

\bibitem{Carlen Loss}
E. Carlen, M. Loss, \textit{Sharp constant in Nash's inequality,}
{International Mathematics Research Notices.} \textbf{7} (1993),
213-215.

\bibitem{Carleson Chang}
L. Carleson, S.Y.A. Chang, \textit{On the existence of an extremal
function for an inequality of J. Moser,} {Bull. Sci. Math.}
\textbf{110} (1986), 113-127.

\bibitem{W Chen}
W. Chen, \textit{$L_{p}$ Minkowski problem with not necessarily
positive data,} {Adv. Math.} \textbf{201} (2006) 77-89.

\bibitem{KS Chou XJ Wang}
 K.S. Chou, X.J. Wang,  \textit{The $L_{
p}-$Minkowski problem and the Minkowski problem in centroaffine
geometry,} {Adv. Math.}  \textit{205} (2006) 33-83.


\bibitem{A Cianchi E Lutwak D Yang G Zhang}
A. Cianchi, E. Lutwak, D. Yang and G. Zhang, \textit{Affine
Moser-Trudinger and Morrey-Sobolev inequalities,} {Calc. Var.
Partial Differ. Equ.} {Doi 10.1007/s00526-009-0235-4}.

 \bibitem{Cordero-Erausquin nazaret Villani}
D. Cordero-Erausquin, B. Nazaret and C. Villani,  \textit{A mass
transportation approach to sharp Sobolev and Gagliardo¨CNirenberg
inequalities,} {Adv. Math.} \textbf{182} (2004), 307-332.


\bibitem{Del-Pino Dolbeaut}
  M. Del Pino,  J. Dolbeault, \textit{Best constants for Gagliardo-Nirenberg
inequalities and applications to nonlinear diffusions,} {J. Math.
Pures Appl.} \textbf{81} (9) (2002), 847-875.


\bibitem{Del-pino Dolbeault 2}
 M. Del Pino, J. Dolbeault, \textit{Nonlinear diffusion and optimal
constants in Sobolev type inequalities: Asymptotic behaviour of
equations involving p-Laplacian,} {C. R. Math. Acad. Sci. Paris,
Ser. I} \textbf{334} (2002),  365-370.


\bibitem{del-pino Dolbeault 1}
 M. Del Pino,  J. Dolbeault, \textit{The optimal Euclidean $L^{p}-$Sobolev
logarithmic inequaity,} {J. Funct. Anal.} \textbf{197} (1) (2003),
151-161.

\bibitem{Del-pino Dolbeault 3}
 M. Del Pino, J. Dolbeault, \textit{Asymptotic behaviour of nonlinear
diffusions,} {Math. Res. Lett.} \textbf{10} (4) (2003),  551-557.

\bibitem{Del-pino Dolbeault Gentil}
 M. Del Pino, J. Dolbeault  and I. Gentil, \textit{Nonlinear diffusions,
hypercontractivity and the optimal $L_{p}-$Euclidean logarithmic
Sobolev inequality,} {J. Math. Anal. Appl.} \textbf{293} (2) (2004)
375-388.


\bibitem{Federer Fleming}
H. Federer, W. Fleming, \textit{Normal and integral currents,} {Ann.
Math.} \textbf{72}  (1960) 458-520.


\bibitem{Gagliardo}
E. Gagliardo, \textit{Propriet\`{a} di alcune classi di funzioni
pi\`{u} variabili,} {Ric. Mat.} \textbf{7} (1958), 102-137.


\bibitem{Gardner}
R.J. Gardner, \textit{Geometric Tomography,} {Encyclopedia of
mathematics and Its Applications, vol. 58,} {Cambridge University
Press,} Cambridge, 1995.

 \bibitem{Gentil} I. Gentil, \textit{The
general optimal $L^{p}$-Euclidean logarithmic Sobolev inequality by
Hamilton-Jacobi equations,} {J. Funct. Anal.} \textbf{202} (2003),
591-599.

\bibitem{Gross}
L. Gross, \textit{Logarithmic Sobolev inequalities,} {Amer. J.
Math.} \textbf{97} (1975), 1061-1083.

\bibitem{K Hilden}
K. Hilden, \textit{Symmetrization of functions in Sobolev spaces and
the isoperimetric inequality,} {Manus. Math.} \textbf{18} (1976),
215-235.

\bibitem{C Hu XN Ma C Shen}
C. Hu, X.N. Ma,C. Shen, \textit{On the Christoffel-Minkowski problem
of Firey's p-sum,} {Calc. Var. Partial Differ. Equ.} \textbf{21}
(2004), 137-155.

\bibitem{B Kawohl1}
B. Kawohl, \textit{Rearrangements and convexity of level sets in
PDE,} {Lecture Notes in Math., vol 1150.} {Springer, Berlin} (1985).


\bibitem{B Kawohl}
B. Kawohl, \textit{On the isoperimetric nature of a rearrangement
inequality and its consequences for some variational problems,}
{Arch. Ration. Mech. Anal.} \textbf{94} (1986), 227-243.



\bibitem{M Ludwig}
M. Ludwig, \textit{Ellipsoids and matrix-valued valuations,} {Duke
Math. J.} \textbf{119}   (2003), 159-188.

\bibitem{M Ludwig 2}
M. Ludwig, \textit{Minkowski valuations.} {Trans. Am. Math. Soc,}
\textbf{357}   (2005) 4191-4213.

\bibitem{E Lutwak 1}
E. Lutwak, \textit{The Brunn-Minkowski-Firey theory. I. Mixed
volumes and the Minkowski problem,} {J. Differ. Geom.} \textbf{38}
(1993), 131-150.


\bibitem{E. Lutwak}
E.  Lutwak, \textit{The Brunn¨CMinkowski-Firey theory. II. Affine
and geominimal surface areas,} {Adv. Math.} \textbf{118} (1996),
244-294.

\bibitem{E Lutwak V Oliker}
E. Lutwak, V. Oliker, \textit{On the regularity of solutions to a
generalization of the Minkowski problem,} {J. Differ. Geom.}
\textbf{41} (1995), 227-246.


\bibitem{E Lutwak Yang Zhang}
E. Lutwak, D. Yang, G. Zhang, \textit{A new ellipsoid associated
with convex bodies,} {Duke Math. J.} \textbf{104}  (2000), 375-390.

\bibitem{Lutwak Yang Zhang}
E. Lutwak, D. Yang, G. Zhang, \textit{On $L^{p}$ affine
isoperimetric inequalities,} {J. Differ. Geom.} \textbf{56} (2000),
111-132.


\bibitem{Lutwak Yang Zhang 1}
E. Lutwak, E. Yang and G. Zhang, \textit{Sharp affine $L_{p}$
Sobolev inequalities,} {J. Differ. Geom.} \textbf{62} (2002), 17-38.


\bibitem{Lutwak yang Zhang 2}
E. Lutwak, D. Yang, G. Zhang, \textit{The Cramer-Rao inequality for
star bodies,} {Duke Math. J.} \textbf{112}  (2002) 59-81.

\bibitem{Lutwak yang Zhang 3}
E. Lutwak, D. Yang, G. Zhang, \textit{On the $L_{p}$ Minkowski
problem,} {Trans. Am. Math. Soc.} \textbf{356},  (2004) 4359-4370.

\bibitem{Lutwak yang Zhang 4}
E. Lutwak, D. Yang, G. Zhang, \textit{Optimal Sobolev norms and the
$L_{p}$ Minkowski problem,} {Int. Math. Res. Not.} (2006) Art. ID
62987, 21 pp.

\bibitem{Mazya}
V.M. Maz'ya, \textit{Classes of domains and imbedding theorems for
function spaces,} {Dokl. Akad. Nauk. SSSR} 133 (1960), 527-530
{(Russian); English translation: Soviet Math. Dokl. 1} (1960),
882-885.

\bibitem{Moser}
J. Moser, \textit{A sharp form of an inequality by N. Trudinger,}
{Indiana. Unv. Math. J.} \textbf{20} (1971), 1077-1092.

\bibitem{Nirenberg}
L. Nirenberg, \textit{On elliptic partial differential equations,}
{Ann. Sc. Norm. Pisa} \textbf{13} (1959), 116-162.

\bibitem{Schneider}
R. Schneider, \textit{Convex Bodies: The Brunn-Minkowski Theory,}
{Encyclopedia of Mathematics and Its Applications, vol 44,}
{Cambridge University Press,} Cambridge, 1993.

\bibitem{E Sperner}
E. Sperner, \textit{Symmetrisierung f\"{u}r Funktionen mehrerer
reeller Variablen,} {Manus. Math.} \textbf{11} (1974), 159-170.

\bibitem{G Talenti}
G. Talenti, \textit{Best constant in Sobolev inequality,} {Ann. Mat.
Pura Appl.} \textbf{110} (1976), 353-372.


\bibitem{A.C. Thompson}
A.C. Thompson,
 \textit{Minkowski Geometry,} {Encyclopedia of
Mathematics and Its Applications, vol. 63,} Cambridge University
Press, Cambridge, 1996.


\bibitem{Weissler}
F.B. Weissler, \textit{Logarithmic Sobolev inequalities for the
heat-diffusion semigroup,} {Trans. Amer. Math. Soc. } \textbf{237}
(1978), 255-269.

\bibitem{G Zhang}
G. Zhang, \textit{The affine Sobolev inequality,} {J. Differential
Geom.} \textbf{53} (1) (1999), 183-202.
\end{thebibliography}
\end{document}